\documentclass[12pt]{article}
\usepackage{euscript,amsmath, amssymb, amsfonts}
\usepackage{color}

\newcommand{\R}{{\mathbb R}}
\newcommand{ \bee}{\begin{eqnarray}}
\newcommand{ \eee}{\end{eqnarray}}
\newcommand\nn{\nonumber \\}
\newcommand{\p}{{\hbar^2}}
\newcommand{\pp}[2]{\hbar^{#1}}

\newcommand{\oC}{{\mathbb C}}

\newcommand{\di}{{ d }}

\newcommand{\G}{{\mathbb G}}
\newcommand{\Z}{{\mathbb Z}}
\newcommand{\K}{{\mathbb K}}
\newcommand{\be}{\begin{equation}}
\newcommand{\ee}{\end{equation}}

\newcommand{\oD}{{\mathbf D}}
\newcommand{\oE}{{\mathbf E}}
\newcommand{\oS}{{\mathbf S}}
\newcommand{\oZ}{{\mathbf Z}}
\newcommand{\eE}{\EuScript E}

\newcounter{sectionn}
\newcommand{\sectionn}[1]{\par\noindent\refstepcounter{sectionn}
                             {\bf \arabic{sectionn}. #1.}}

\newcounter{theorem}
\newcommand{\theorem}{\par\refstepcounter{theorem}
           {\bf Theorem \arabic{section}.\arabic{theorem}. }}

\makeatletter \@addtoreset{theorem}{section}

\makeatletter \@addtoreset{equation}{section}

\makeatother

\font\frtnfr=eufm10   scaled\magstep1
\font\twlfr=eufm10
\font\tenfr=eufm10

\newfam\frfam
\textfont\frfam=\frtnfr
\scriptfont\frfam=\twlfr
\scriptscriptfont\frfam=\tenfr

\font\frtnopen=msbm10  scaled\magstep2
\font\twlopen=msbm10
\font\tenopen=msbm10

\newfam\openfam
\textfont\openfam=\frtnopen
\scriptfont\openfam=\twlopen
\scriptscriptfont\openfam=\tenopen
\def\open{\fam\openfam}

\font\frtnsf = cmss12 scaled\magstep1
\font\twlsf = cmss10
\font\tensf = cmss9

\newfam\Scfam
\textfont\Scfam = \frtnsf
\scriptfont\Scfam = \twlsf
\scriptscriptfont\Scfam = \tensf

\begin{document}

%\mark{ DCQIS-07 \ddatt}

\sloppy \title
 {
The deformations of nondegenerate constant Poisson bracket
with even and odd deformation parameters}
\author
 {
 S.E.Konstein\thanks{E-mail: konstein@lpi.ru}\ \ and
 I.V.Tyutin\thanks{E-mail: tyutin@lpi.ru}
 \thanks{
               This work was supported
               by the RFBR (grant No.~08-02-01118),
               and by the grant LSS-1615.2008.2
 } \\
               {\sf \small I.E.Tamm Department of
               Theoretical Physics,} \\ {\sf \small P. N. Lebedev Physical
               Institute,} \\ {\sf \small 119991, Leninsky Prospect 53,
               Moscow, Russia.} }
\date {
%\ddatt
}

\maketitle

\begin{abstract}
{ \footnotesize We consider Poisson superalgebras with constant nondegenerate
bracket realized on the smooth
Grassmann-valued functions with compact supports in $\R^n$. The
 deformations with even and odd deformation parameter
 of these superalgebras
are presented for $n\ge 4$.  }
\end{abstract}

%%%%%%%%%%%%%%%%%%%%%%%%%%%%%%%%%%%%%%%%%%%%%%%%%%%

\section{Introduction.}

The hope to construct the quantum mechanics on nontrivial
manifolds is connected with geometrical or deformation
quantization \cite{1} - \cite{4}.
The functions on the phase space
are associated with the operators, and the product and the
commutator of the operators are described by associative *-product
and *-commutator of the functions. These *-product and
*-commutator are the deformations of usual product and of
usual Poisson bracket. The deformations of Poisson (anti)bracket
was considered in many publications for different spaces of functions.
In \cite{LS}, the problem is considered for the
Poisson superalgebra on the superspace of polynomials and for the antibracket.
Purely Grassmannian case is considered in \cite{Ty2}.
Bosonic Poisson algebra realized on smooth functions was considered in \cite{Zh}.

The result depends on chosen function space (see \cite{LS} for more details).

It occurred, that there exist odd second cohomologies with
coefficients in adjoint representation. It is natural to look for
the deformations associated with these odd cohomologies and having
the odd deformation parameter \cite{LS}, \cite{leit}.

\sectionn{Deformations of topological Lie superalgebras}
We
recall some concepts concerning formal deformations of algebras (see,
e.g.,~\cite{Gerstenhaber}), adapting them to the case of topological Lie
superalgebras.
Let $L$ be a topological Lie superalgebra over $\K$ ($\K=\R$
or $\oC$) with Lie superbracket $\{\cdot,\cdot\}$, $\K[[\p]]$ be the ring of
formal power series in $\p$ over $\K$, and $L[[\p]]$ be the $\K[[\p]]$-module
of formal power series in $\p$ with coefficients in $L$.
We endow both
$\K[[\p]]$ and $L[[\p]]$ by the direct-product topology.
The grading of $L$
naturally determines a grading of $L[[\p]]$:  an element $f=f_0+\p
f_1+\ldots$ has a definite parity $\varepsilon(f)$ if
$\varepsilon(f)=\varepsilon(f_j)$ for all $j=0,1,..$.
Every $p$-linear
separately continuous mapping from $L^p$ to $L$ (in particular, the bracket
$\{\cdot,\cdot\}$) is uniquely extended by $\K[[\p]]$-linearity to a
$p$-linear separately continuous mapping over $\K[[\p]]$ from $L[[\p]]^p$ to
$L[[\p]]$.
A (continuous) formal deformation of $L$ is by definition a
$\K[[\p]]$-bilinear separately continuous Lie superbracket $C(\cdot,\cdot)$
on $L[[\p]]$ such that $C(f,g)=\{f,g\} \mod \p$ for any $f,g\in L[[\p]]$.
Obviously, every formal deformation $C$ is expressible in the form
\begin{equation} C(f,g)=\{f,g\}+\p
C_1(f,g)+\pp{4}{2}C_2(f,g)+\ldots,\quad f,g\in L,
\label{1}
\end{equation}
where $C_j$
are separately continuous skew-symmetric bilinear mappings from $L\times L$
to $L$ (2-cochains with coefficients in the adjoint representation of $L$).
Formal deformations $C^1$ and $C^2$ are called equivalent if there is a
continuous $\K[[\p]]$-linear operator $T: L[[\p]]\to L[[\p]]$ such that
$TC^1(f,g)=C^2(T f,Tg)$, $f,g\in L[[\p]]$ and $T=\mbox{ id}+\p T_1$.
  The problem of finding formal
deformations of $L$ is closely related to the problem of computing
Chevalle--Eilenberg cohomology of $L$ with coefficients in the adjoint
representation of $L$.  Let $\mathcal C_p(L)$ denote the space of $p$-linear
skew-symmetric separately continuous mappings from $L^p$ to $L$ (the space of
$p$-cochains with coefficients in the adjoint representation of $L$).  The
space $\mathcal C_p(L)$ possesses a natural $\Z_2$-grading: by definition,
$M_p\in \mathcal C_p(L)$ has the definite parity $\epsilon(M_p)$ if the
relation
$$
\varepsilon(M_p(f_1,\ldots,f_p))=
\varepsilon(M_p)+\varepsilon(f_1)+\ldots+\varepsilon(f_1)
$$
holds for any $f_j\in L$ with definite parities $\epsilon(f_j)$.
Since a Lie superbracket is always even,
all $C_j$ in the expansion~(\ref{1}) should be even 2-cochains.
The differential $\di_p^{\rm ad}$ is defined to be the linear
operator from $\mathcal C_p(L)$ to $\mathcal C_{p+1}(L)$ such that
\begin{align}
&&d_p^{\rm ad}M_p(f_1,...,f_{p+1})=
-\sum_{j=1}^{p+1}(-1)^{j+\varepsilon(f_j)|\varepsilon(f)|_{1,j-1}+
\varepsilon(f_j)\varepsilon_{M_p}}\{f_j,
M_p(f_{1},...,\hat{f}_j,...,f_{p+1})\}- \nonumber \\
&&-\sum_{i<j}(-1)^{j+\varepsilon(f_j)|\varepsilon(f)|_{i+1,j-1}}
M_p(f_1,...f_{i-1},\{f_i,f_j\},f_{i+1},...,\hat{f}_j,...,f_{p+1}),\nonumber
\end{align}
for any $M_p\in \mathcal C_p(L)$ and $f_1,\ldots f_{p+1}\in L$
having definite parities.
Here the hat means that the argument is omitted and the notation
$
|\varepsilon(f)|_{i,j}=\sum_{l=i}^j\varepsilon(f_l)
$
has been used.
Writing the Jacobi identity
for a deformation $C$ of the form (\ref{1}),
\begin{equation}
\nonumber
(-1)^{\varepsilon(f)\varepsilon(h)}C(f,C(g,h))+\mathrm{cycle}(f,g,h)=0,
\end{equation}
and taking the terms of the order $\p$, we find that
\begin{equation}
\nonumber
d_2^{\rm ad}C_1=0.
\end{equation}
Thus, the
first order deformations of $L$
are described by 2-cocycles of the differential $d^{\rm ad}$.

\sectionn{Poisson superalgebra}
Let $\EuScript D(\R^k)$ denote the space of
smooth $\K$-valued functions with compact support on $\R^k$.
This space is
endowed by its standard topology: by definition, a sequence $\varphi_k\in
\EuScript D(\R^k)$ converges to $\varphi\in \EuScript D(\R^k)$ if
$\partial^\lambda\varphi_k$ converge uniformly to $\partial^\lambda\varphi$
for every multi-index $\lambda$, and the supports of all $\varphi_k$ are
contained in a fixed compact set. We set
$$  D^{n_-}_{n_+}= \EuScript
D(\R^{n_+})\otimes \G^{n_-},\quad  E^{n_-}_{n_+}=
C^\infty(\R^{n_+})\otimes \G^{n_-},
$$
where $\G^{n_-}$ is the Grassmann
algebra with $n_-$ generators.  The generators of
the Grassmann algebra (resp., the coordinates of the space $\R^{n_+}$) are
denoted by $\xi_\alpha$, $\alpha=1,\ldots,n_-$ (resp., $x_i$, $i=1,\ldots,
n_+$).  We shall also use collective variables $z_A$ which are equal to $x_A$
for $A=1,\ldots,n_+$ and are equal to $\xi_{A-n_+}$ for
$A=n_++1,\ldots,n_++n_-$.  The spaces $ D^{n_-}_{n_+}$ and
$E^{n_-}_{n_+}$ possess a natural grading
which is determined by that of the Grassmann algebra. The parity of an
element $f$ of these spaces is denoted by $\varepsilon(f)$. We also set
$\varepsilon_A=0$ for $A=1,\ldots, n_+$ and $\varepsilon_A=1$ for
$A=n_++1,\ldots, n_++n_-$.

Let $\partial/\partial z_A$ and $\overleftarrow{\partial}/\partial z_A$ be
the operators of the left and right differentiation.  The Poisson bracket is
defined by the relation
\begin{equation}
\{f,g\}(z)=\sum_{A,B} f(z)\frac{\overleftarrow{\partial}}{\partial z_A}\omega^{AB}
\frac{\partial}{\partial z_B}g(z)=
-(-1)^{\varepsilon(f)\varepsilon(g)}\{g,f\}(z),\label{3.0}
\end{equation}
where
the symplectic metric $\omega^{AB}=(-1)^{\varepsilon_A
\varepsilon_B}\omega^{BA}$ is a constant invertible matrix.  For
 definiteness, we choose it in the form
\[
\omega^{AB}=
\left(\begin{array}{cc}\omega^{ij}&0       \\
0&\lambda_\alpha\delta^{\alpha\beta}\end{array}\right),\quad
\lambda_\alpha=\pm1,\ i,j=1,...,n_+,\ \alpha,\beta=1,...,n_-
\]
 where
$\omega^{ij}$ is the canonical symplectic form (if $\K=\oC$, then one can
choose $\lambda_\alpha=1$).  The nondegeneracy of the matrix $\omega^{AB}$
implies, in particular, that $n_+$ is even.  The Poisson superbracket
satisfies the Jacobi identity
\begin{equation}
(-1)^{\varepsilon(f)\varepsilon(h)}
\{f,\{g,h\}\}(z)+\hbox{cycle}(f,g,h)= 0,\quad f,g,h\in E^{n_-}_{n_+}.
\label{3.0a}
\end{equation}
By Poisson superalgebra $\mathcal
P^{n_-}_{n_+}$, we mean the space $D^{n_-}_{n_+}$ with the Poisson
bracket~(\ref{3.0}) on it.  The relations~(\ref{3.0}) and~(\ref{3.0a}) show
that this bracket indeed determines a Lie superalgebra structure on $
D^{n_-}_{n_+}$.

The integral on $ D^{n_-}_{n_+}$ is defined by the relation
$$
\bar
f\stackrel{\mathrm{def}}{=}\int \di z\, f(z)= \int_{\R^{n_+}}\di x\int
\di\xi\, f(z),
$$
where the integral on the Grassmann algebra is normed by
the condition $\int \di\xi\, \xi_1\ldots\xi_{n_-}=1$.

Introduce the superalgebra $Z^{n_-}_{n_+}$,
$D^{n_-}_{n_+}\subset Z^{n_-}_{n_+}\subset E^{n_-}_{n_+}$,
$Z^{n_-}_{n_+}=D^{n_-}_{n_+}\oplus {\cal C}_{D^{n_-}_{n_+}}(E^{n_-}_{n_+})$,
where ${\cal C}_{D^{n_-}_{n_+}}(E^{n_-}_{n_+})$ is centralizer of
$D^{n_-}_{n_+}$ in $E^{n_-}_{n_+}$.

Introduce the notation
\be
{ S}_n^m\stackrel {def} =
 { E}_n^m/{ Z}_n^m.
 \ee

\sectionn{Poisson superalgebra with outer odd parameters}
Below we consider the case, where the functions and multilinear forms may
depend on outer odd parameters $\theta_i$, where $\theta$-s belong to some supercommutative
associative superalgebra. For simplicity we consider the case $\theta_i\in {\open G}^k$.
Thus we consider a colored algebra
$\oD=\G^{k}\otimes D^{n_-}_{n_+}$, with $(\Z_2)^2$ grading, namely
the grading of element $\theta\otimes f$ is $(\varepsilon_1 (\theta),
\varepsilon_2 (f))$.

Below we consider $\oD$ as a Lie superalgebra with the parity $\varepsilon=\varepsilon_1
+\varepsilon_2$. One can easily check that such consideration is selfconsistent
(see also \cite{Sche1} and discussion on Necludova and Sheunert theorems in \cite{Sosrus1}.

Analogously we introduce superspaces and superalgebras
\bee
\oE=\G^k\otimes E_{n_+}^{n_-} \nonumber \\
\oZ=\oD\oplus {\cal C}_{\oD}(\oE) \nonumber \\
\oS=\G^k\otimes S_{n_+}^{n_-} \nonumber \\
\eee

In the following sections we will consider two cases
separately:

\begin{enumerate}
\item $k=1$, there exists only one odd parameter $\theta$, $\theta^2=0$.
\item $k>1$.
\end{enumerate}

All these parameters $\theta$ are supposed to be generating elements of $\G^k$,
i.e. they are odd and satisfy the following condition:

if $\theta y =0$ for some $y\in \G^k$ then $y=\theta z$ for some $z\in \G^k$.

\sectionn{Jacobiators}
Let $p$, $q$ be antisymmetric bilinear forms.
Here and below Jacobiators are defined as follows:
\begin{eqnarray*}
J(p,q)&\stackrel {def} = &
(-1)^{\varepsilon(f)\varepsilon(h)}
\left (
p(q(f,g), h)+q(p(f,g),h)\right) +\mbox{cycle}(f,g,h),
\\
J(p,p)&\stackrel {def} = &
(-1)^{\varepsilon(f)\varepsilon(h)}p(p(f,g),h)+\mbox{cycle}(f,g,h).
\end{eqnarray*}
Evidently, $J(p,q)\in {\cal C}_3(\oD , \oD)$,
if $p,q\in {\cal C}_2(\oD , \oD)$.

If $m_0(f,g)=\{f,g\}$ then $-(-1)^{\varepsilon(f)\varepsilon(h)}J(p,m_0)=d_2^{\mathrm {ad}}p$.

\sectionn{Sign rules}
We use here the following sign rules for factorization the odd parameters
\be
M_n(\theta f_1, f_2,...,f_n)
=(-1)^{\varepsilon (\theta)\varepsilon(M_n)}\theta M_n(f_1,..., f_n)
\ee
It follows from this sign rule and superantisymmetry of the form $M_n$
that`
\bee
M_n(f_1, \theta f_2, ... ,f_n)= M_n (f_1 \theta, f_2, ... , f_n)\\
M_n(f_1, f_2, ... , f_n \theta)= M_n (f_1, f_2, ... , f_n)\theta.
\eee

\section{1 odd deformation parameter}

Let $\theta$ be the odd deformation parameter, $\varepsilon (\theta)=1$.

Let deformation $C(f,g)$ has the form
$$C(f,g) = C_0(f,g)+\theta C_1(f,g),$$
where $\varepsilon (C_0)=0$, $\varepsilon (C_1)=1$.

Jacoby identity gives
\be
0=J(C,C)=J(C_0, C_0) +J(C_0, \theta C_1),
\ee
which implies
\bee
J(C_0, C_0)=0\nn
J(C_0, \theta C_1)=0 \label{main}
\eee

So, $C_0$ is a deformation of the Poisson superalgebra,
and (\ref{main}) is the equation, which we will investigate
below for some of the deformations $C_0$, found in \cite{anti} and \cite{deform4}.

\subsection{Antibracket}

An interesting example of Lie superalgebra with 1 even and 1 odd cohomology
is antibracket realized on $D_n^n$.

The spaces $ D^{n}_{n}$ and
$E^{n}_{n}$ possess also another $\Z_2$-grading
$\epsilon$ ($\epsilon$-parity),
 which is inverse to $\varepsilon$-parity: $\epsilon=\varepsilon+1$.

We set
$\varepsilon_A=0$, $\epsilon_A=1$ for $A=1,\ldots, n_+$
and $\varepsilon_A=1$, $\epsilon_A=0$ for
$A=n_++1,\ldots, n_++n_-$.

It is well known, that the bracket
\bee
\nonumber
[f,g](z)=\sum_{i=1}^n\left(f(z)\frac{\overleftarrow{\partial}}{\partial x_i}
\frac\partial{\partial\xi_i}g(z)-
f(z)\frac{\overleftarrow{\partial}}{\partial\xi_i}
\frac{\partial}{\partial x_i}g(z)\right),
\eee
which we following to \cite{BV1} call "antibracket" or "odd bracket",
defines the structure of Lie superalgebra on the superspaces
$D_{n}\stackrel {def} =  D^{n}_{n}$
and $ E_{n}\stackrel {def} =  E^{n}_{n}$
with the $\epsilon$-parity.

Indeed, $[f,g]=-(-1)^{\epsilon(f)\epsilon(g)}[g,f]$,
$\epsilon([f,g])=\epsilon(f)+\epsilon(g)$,
and Jacobi identity is satisfied:
\be
\nonumber
(-1)^{\epsilon(f)\epsilon(h)}[f,[g,h]]
+(-1)^{\epsilon(g)\epsilon(f)}[g,[h,f]]
+(-1)^{\epsilon(h)\epsilon(g)}[h,[f,g]]=0
,\quad f,g,h\in E_{n}.
\ee

Here these Lie superalgebras are called antiPoisson
superalgebras.

The odd Poisson bracket play an important role in Lagrangian formulation of
the quantum theory of the gauge fields, which is known as BV-formalism
\cite{BV1}, \cite{BV2} (see also \cite{reports}-\cite{HenTei}, \cite{LS}).
These odd
bracket were introduced in physical literature in \cite{BV1}.
Antibracket possesses many features analogous to ones
of even Poisson bracket and even can be obtained via "canonical formalism"
with odd time. However, contrary to the case of even Poisson bracket where
there exists voluminous literature on different aspects of the deformation
(quantization) of Poisson algebra, the problem of the deformation of
antibracket is not study satisfactory yet.

\theorem\label{thdef}{}\cite{anti}
{\it Up to similarity transformation,
the deformation of antiPoisson superalgebra with even parameter $\hbar$
has the form
\begin{equation}
\label{def}
[f(z),\, g(z)]_*\!=\![f(z),g(z)] + (-1)^{\varepsilon(f)}\{\frac{
c}{1+ cN_z/2}\Delta f(z)\} \eE_zg(z)
+
\{\eE_z f(z)\}
\frac{ c}{1+ cN_z/2}\Delta g(z),
\end{equation}
where $N_z=\sum_A z_A \frac {\partial}{\partial z_A}$,
$\Delta
=\sum_i\frac \partial {\partial x_i}\frac \partial {\partial \xi_i}
$
and $c\in \p\K[[\p]]$.
}

\subsubsection{Odd cohomology of antiPoisson superalgebra and corresponding deformation}

Odd cohomology of antiPoisson superalgebra has the form \cite{anticoh}
\be
m_{2|3}(z|f,g)=(-1)^{\varepsilon (f)}\{(1-N_{\xi })f(z)\}(1-N_{\xi })g(z)
\ee

One can prove that in this case the  equation (\ref{main})
has the only solution which leads to the
following deformation
$$[f,g]_\star = [f,g]+\theta \left(1-\xi^i \partial_{\xi^i}
\right) f\cdot \left(1-\xi^i \partial_{\xi^i} \right)g$$

\subsection{Deformations of the Poisson superalgebra.}

For any $\varkappa\in \K[[\hbar]]$, such that $c_1\stackrel {def}=\frac 1 6
\p\varkappa^2\in
\p\K[[\p]]$, the
Moyal-type superbracket
\begin{equation}
\nonumber
{\cal M}_{c_1}(z|f,g)=\frac{1}{\hbar\varkappa}f(z)\sinh
\left(\hbar\varkappa\sum_{A,B}\frac{\overleftarrow{\partial}}{\partial z_A}\omega^{AB}
\frac{\partial}{\partial z_B}\right)g(z)
\end{equation}
is skew-symmetric and
satisfies the Jacobi identity
and, therefore, gives a deformation of the
initial Poisson algebra.

Let the bilinear mappings
$m_3$,
 and $m_\zeta$ from
 $(D^{n_-}_{n_+})^2$ to $ D^{n_-}_{n_+}$ be defined by the relations
\begin{eqnarray}
m_3(z|f,g)& = & (-1)^{n_-\varepsilon(f)}\eE_z f(z)\bar g -
(-1)^{\varepsilon(f)\varepsilon(g)+n_-\varepsilon(g)}  \eE_z g(z)\bar f,       \label{m3}
\\
\label{mzeta}
m_\zeta(z|f,g)&=&(-1)^{n_-\varepsilon(f)} \{\zeta(z),f(z)\} \bar{g} -
  (-1)^{\varepsilon(f)\varepsilon(g) + n_-\varepsilon(g)}\{\zeta(z),g(z)\} \bar{f},
\end{eqnarray}
where
$
\eE_z  \stackrel {def} =  1-\frac 1 2 z \partial_z,    $
and
$z = (x_1,...,x_n,\xi_1,...,\xi_m)$.

For $\zeta\in E^{n_-}_{n_+}[[\p]]$, $c_1, c\in \K[[\p]]$,
we set
\begin{align} &C^{(1)}_{\zeta, c_1}(z|f,g)=
{\cal M}_{c_1}(z|f+\zeta\bar{f},g+\zeta\bar{g}),
\nonumber\\
&C^{(1)}_{\zeta,c_1,c}(z|f,g)= {\cal
M}_{c_1}(z|f+\zeta\bar{f},g+\zeta\bar{g})+c\bar{f}\bar{g}\nonumber
\end{align}

\subsubsection{Deformations of the Poisson superalgebra at $ n_+\ge 4$.}

\medskip

{\begin{theorem}\label{theo}
\cite{deform4}
\it Let $\theta_i=0$ for all $i$. Then
\begin{enumerate}
\item Let $n_-=2k$ and $ n_+\ge 4$. Then every
continuous formal deformation of the Poisson superalgebra $D^{n_-}_{n_+} $
is
equivalent either to the superbracket $C^{(1)}_{\zeta,c_1}(z|f,g)$, where
$\zeta\in \p E^{n_-}_{n_+}[[\p]]$ is even and $c_1\in \p\K[[\p]]$, or
to the superbracket
$$
C^{(3)}_{\zeta, c_3}(z|f,g)=\{f(z),g(z)\}+m_{\zeta}(z|f,g)+ c_3 m_3(z|f,g),
$$
where $\zeta\in \p E^{n_-}_{n_+}[[\p]]$ is even and $c_3\in \p\K[[\p]]$.
The deformations $C^{(i)}_{\zeta_1, c}$ and $C^{(i)}_{\zeta_2, c}$ are
equivalent if $\zeta_1 - \zeta_2 \in Z^{n_-}_{n_+}[[\p]]$.

\item Let $n_-=2k+1$ and $n_+\ge 4$. Then every continuous formal deformation of the Poisson
superalgebra $D^{n_-}_{n_+}$ is equivalent to the superbracket
$C^{(1)}_{\zeta,c_1,c}(z|f,g)$, where $c,c_1\in \p\K[[\p]]$ and $\zeta\in
\p E^{n_-}_{n_+}[[\p]]$ is an odd function such that $[{\cal
M}_{c_1}(z|\zeta,\zeta)+ c]\in  D^{n_-}_{n_+}[[\p]]$.
The deformations $C^{(1)}_{\zeta_1, c_1,c}$ and $C^{(1)}_{\zeta_2, c_1, c}$ are
equivalent if $\zeta_1 - \zeta_2 \in Z^{n_-}_{n_+}[[\p]]$.

\end{enumerate}
\end{theorem}}

\subsubsection{Odd cohomology}

\begin{enumerate}
\item Let $n_-=2k$ and $ n_+\ge 4$.

Then the odd cohomology has the form
\be
m_{\zeta}(z|f,g)
\ee
with odd function $\zeta$.

\item Let $n_-=2k+1$ and $n_+\ge 4$.

Then the odd cohomology has either the form
\be
m_{\zeta}(z|f,g)
\ee
with even function $\zeta$
or the form
\be
m_3(z|f,g).
\ee
\end{enumerate}

\subsubsection{Deformation of Poisson superalgebra for even $n_-$}

\medskip

{\begin{theorem} \it
Let $n_-$ be even and $n_+\ge 4$. Then each deformation of Poisson
superalgebra with one odd deformation parameter is equivalent either to
$$
\{f,g\}_\star\,=\,
\frac{1}{\hbar\kappa}(f+\zeta \bar f)\sinh
\left(\hbar\kappa\frac{\overleftarrow{\partial}}{\partial z^A}\omega^{AB}
\frac{\partial}{\partial z^B}\right)(g+\zeta \bar g)
$$
or to
$$
\{f,g\}_\star=\,\{f,g\}+m_\zeta(f,g)+cm_3(f,g),
$$
where
$\zeta =\zeta_0+\theta \zeta_1,\ \zeta_i \in\hbar^2
{S}_{n_+}^{n_-}[[\hbar^2]],\ \epsilon(\zeta_i)=i$
and
$\kappa^2,\, c \in \hbar^2{\mathbf C}[[\hbar^2]]$.
\end{theorem}
}

This Theorem can be proved analogously to Theorem \ref{theo}.

\subsubsection{Deformation of Poisson superalgebra for odd $n_-$}

\medskip

{\begin{theorem} \it
Let $n_-$ be odd and $n_+\ge 4$. Then each deformation of Poisson algebra with one
odd deformation parameter is equivalent either to
$$
\{f,g\}_\star=
\frac{1}{\hbar\kappa}(f+\zeta \bar f)\sinh
\left(\hbar\kappa\frac{\overleftarrow{\partial}}{\partial z^A}\omega^{AB}
\frac{\partial}{\partial z^B}\right)(g+\zeta \bar g)+c\,\bar f\,\bar g
$$
or to
$$
\{f,g\}_\star=\,\{f,g\}+
m_\zeta(f,g)+\theta m_3(f,g),$$
where
$$\zeta=\zeta_1+\theta \zeta_0,\
\zeta_i\in \hbar^2{S}_{n_+}^{n_-}[[\hbar^2]], \ \epsilon(\zeta_i)=i, $$
$$\kappa^2,c\in \hbar^2{\mathbf C}[[\hbar^2]],$$
are such that
$$
\frac{1}{\hbar\kappa}\zeta\sinh
\left(\hbar\kappa\frac{\overleftarrow{\partial}}{\partial z^A}\omega^{AB}
\frac{\partial}{\partial z^B}\right)\zeta +c \in {\mathbf D}[[\hbar^2]].
$$
\end{theorem}}
This Theorem can be proved analogously to Theorem \ref{theo}.

\section{Finite number of odd deformation parameters}

Let us look for deformation of Poisson superalgebra with more than 1
odd deformation parameters.

\subsection{Deformation of Poisson superalgebra for even $n_-$}

In the case of even $n_-$ we can reduce some number of odd cohomologies
$m_{\zeta_{\alpha_1,...,\alpha_k}}$ ($\alpha_i=0,1$)
with odd $\zeta_{\alpha_1,...,\alpha_k}\in S_{n_+}^{n_-}$ to even $m_\zeta$
with even $\zeta\in\oS$:
\be\label{zeta}
\zeta=\sum_s\sum_{{\alpha_1,...,\alpha_k}=0,1}
(\theta_1)^{\alpha_1}\cdot ...\cdot (\theta_k)^{\alpha_k}\cdot
\zeta_{\alpha_1,...,\alpha_k},\ \ \varepsilon(\zeta)=0
\ee

Then one can easily prove the following Theorem, which is complete analog of
the first item of Theorem \ref{theo}.

\medskip

{\begin{theorem}
\it
Let $n_-=2k$ and $ n_+\ge 4$. Then every
continuous formal deformation of the Poisson superalgebra $\oD$
is
equivalent either to the superbracket $C^{(1)}_{\zeta,c_1}(z|f,g)$, where
$\zeta\in \p \oD[[\p]]$ is even and $c_1\in \p\G^k[[\p]]$ is even, or
to the superbracket
$$
C^{(3)}_{\zeta, c_3}(z|f,g)=\{f(z),g(z)\}+m_{\zeta}(z|f,g)+ c_3 m_3(z|f,g),
$$
where $\zeta\in \p \oE[[\p]]$ is even and $c_3\in \p\G^k[[\p]]$ is even.
The deformations $C^{(i)}_{\zeta_1, c}$ and $C^{(i)}_{\zeta_2, c}$ are
equivalent if $\zeta_1 - \zeta_2 \in \oZ[[\p]]$.

\end{theorem}}

\subsection{Deformation of Poisson superalgebra for odd $n_-$}

There are infinite number of odd cohomologies in this case also:
$m_3$ and $m_\zeta$ with even $\zeta\in\S_{n_-}^{n_+}$. We can reduce
some finite part of these odd forms to even by multiplying them by
odd elements of $\G^k$:
$m_3 \to \theta_1 m_3$,
$m_{\zeta_{\alpha_1,...,\alpha_k}\in S_{n_+}^{n_-}} \to m_\zeta$
with odd $\zeta\in\oS$:
\be\label{zeta2}
\zeta=\sum_s\sum_{{\alpha_1,...,\alpha_k}=0,1}
(\theta_1)^{\alpha_1}\cdot ...\cdot (\theta_k)^{\alpha_k}\cdot
\zeta_{\alpha_1,...,\alpha_k},\ \ \varepsilon(\zeta)=1
\ee

Let us look for the deformation. Consider the zero and first order terms:
$$
C(f,g)=m_0(f,g)+h m_1(f,g) + \theta_1 m_3(f,g) + m_\zeta (f,g) +(higher\ order\ terms).
$$

Jacobi identity gives the relation $\theta_1 h=0$.

We restrict ourselves here to the case $\theta_1\ne 0$, which implies $h=\theta_1 h_1$.
Below we omit everywhere the subscript $1$ at $\theta_1$.

The following decomposition orders lead to the next form of the deformation
$$
C=m_{0}+\theta {h}_{1}m_{1}+\theta m_{3}+m_{\zeta }+\theta
h_{1}j_{\zeta }+\eta (z)\mu,
$$
where
\bee
&&j_{\zeta }(f,g)=(-1)^{n_{-}\varepsilon (f)}m_{1}(\zeta ,f)\overline{g}%
-(-1)^{\varepsilon (f)\varepsilon (g)+n_{-}\varepsilon (g)}m_{1}(\zeta ,g)%
\overline{f}         \nonumber\\
&&\mu (f,g)=(-1)^{\varepsilon (f)}\overline{f}\overline{g},\;\varepsilon(\mu)=0, \nonumber
\eee
and even function $\eta(z)\in\oD$ should be determine.

The Jacobi identity gives the following relations
\bee
&&
\{ \zeta, \eta\}+\bar\eta\eta+\theta[\eE  -(2+n_+-n_-)]\eta+h m_1(\zeta,\eta)=0,
\nonumber \\
&&
\theta h=0, \ \ \theta\bar\eta=0, \nonumber \\
&&
h m_1(\{\zeta,\zeta\},f(z))+hm_1(\eta,f(z))+h\bar\eta m_1(\zeta,f(z))=0, \nonumber \\
&&
\eta + hm_1(\zeta,\zeta)+ \theta [2\eE -(2+n_+-n_-)]\zeta+\bar\eta\zeta+\{\zeta,\zeta\}=h_2,
\nonumber
\eee
where $h_2\in\G^k$ is even.

This system of relations is equivalent to the following relations
\bee
&&
\eta=  -\theta h_1m_1(\zeta,\zeta)- \theta [2\eE -(2+n_+-n_-)]\zeta-\bar\eta\zeta-
\{\zeta,\zeta\}+h_2,\nonumber \\
&&
\theta\bar\eta = 0, \nonumber \\
&&
\theta (1+n_+-n_-)h_2 -\bar\eta h_2=0. \nonumber
\eee

Thus, we've obtained theorem

{\it
\begin{theorem}
Let $n_-$ be odd and $n_+\ge 4$. Then Poisson superalgebra
has the deformation depending on $k$ odd parameters
$$
C=m_{0}+\theta {h}_{1}m_{1}+\theta m_{3}+m_{\zeta }+\theta
h_{1}j_{\zeta }+\eta (z)\mu,
$$
where $\zeta\in\oS$, $\eta\in\oD$ and $h_1,\,h_2\in\G^k$ satisfy the relations
\bee
&&
\eta=  -\theta h_1 m_1(\zeta,\zeta)- \theta [2\eE -(2+n_+-n_-)]\zeta-\bar\eta\zeta-
\{\zeta,\zeta\}+h_2,\nonumber \\
&&
\theta\bar\eta = 0, \nonumber \\
&&
\theta (1+n_+-n_-)h_2 -\bar\eta h_2=0, \nonumber\\
&&\varepsilon (z)=1,\ \varepsilon (\eta)=0,\ \varepsilon (h_1)=1,
\ \varepsilon (h_2)=0. \nonumber
\eee
\end{theorem}
}

\vskip 5mm

{\bf Acknowledgements.} Authors are grateful to D.Leites and O.Ogievetski
for useful discussions.


\begin{thebibliography}{99}

\bibitem{1} {\it F.Bayen, M.Flato, C.Fronsdal, A.Lichnerovich,
D.Sternheimer}, Ann.Phys, {\bf 111} (1978) 61; Ann.Phys, {\bf 111} (1978)
111.

\bibitem{2} {\it M.Karasev, V.Maslov}, Nonlinear Poisson brackets. Geometry
and Quantization. Mos\-cow: Na\-uka, 1991.

\bibitem{3} {\it B.Fedosov}, Deformation Quantization and Index Theory.
Berlin: Akademie, 1996.

\bibitem{4} {\it M.Kontsevich}, Deformation Quantization of Poisson
Manifolds, I, q-alg/9709040.

\bibitem{LS}
{\it D.Leites, I.Shchepochkina},
How to quantize the antibracket,
Teor.~Mat.~Fiz., {\bf 126} (2001) 339.

\bibitem{Ty2} {\it I.V.Tyutin}, Teor.~Mat.~Fiz., {\bf 128} (2001) 515.

\bibitem{Zh} {\it V.~V.~Zharinov},  Theor.~Math.~Phys., {\bf 136}, 1049--1065 (2003).

\bibitem{leit} {\it D.Leites}, Clifford algebra as a superalgebra and quantization,
Theor.~Mat.~Fiz., {\bf 58} (1984), no. 2, 229--232.

\bibitem{Sche1} {\it M.~Sheunert}, Generalized Lie Algebras,
J.~Math.~Phys. {\bf 20} (1979),  no.4, 712-720.

\bibitem{Sosrus1} {\it D.Leites} (ed.),
"Seminar on Supersymmetry", vol. 1 "Algebra and Calculus"
(J.Bernstejn, D.Leites, V.Shander, V.Molotkov),
in preparation.

\bibitem{anti} {\it
S.E.Konstein and I.V.Tyutin},
Deformations and central extensions of the antibracket Superalgebra,
Journal of Mathematical Physics, {\bf 49}, 072103 (2008).

\bibitem{deform4} {\it S.~E.~Konstein, A.~G.~Smirnov and I.~V.~Tyutin,}
General form of the deformation of the Poisson superbracket,
Teor.~Mat.~Fiz., {\bf 148} (2006), 1011;
hep-th/0401023.

\bibitem{Gerstenhaber} {\it M.~Gerstenhaber}, Ann. Math. \textbf{79} (1964),
59--103; ibid. \textbf{99} (1974), 257--276.

\bibitem{BV1}
{\it Batalin I.A., Vilkovisky G.A.},
Feynman rules for reducible gauge theories,
Phys. Lett., {\bf 120B} (1983) 166

\bibitem{BV2}
{\it I.A. Batalin and G.A. Vilkovisky},
Existence theorem for gauge algebra,
J. Math. Phys., {\bf 26} (1985) 172

\bibitem{reports}
{\it Gomis J., Paris J., Samuel S.},
Antibrackets, antifields and gauge theory
quantization, Phys. Rep., 1995, {\bf 259}, 1 -- 145 (hep-th/9412228).

\bibitem{books}
{\it D.M.Gitman and I.V.Tyutin},
Quantization of Fields with Constraints,
(Springer--Verlag, 1990).

\bibitem{HenTei}
{\it Henneaux M. and Teitelboim C.},
Quantization of Gauge Systems,
Princeton University Press, Princeton, 1992.

\bibitem{anticoh} {\it S.~E.~Konstein and I.~V.~Tyutin},
Cohomology of antiPoisson superalgebra,
hep-th/0512300.

\end{thebibliography}
\end{document}